\documentclass[12pt,a4paper]{article}
\usepackage{amssymb}

\usepackage{graphicx,amssymb,amsfonts,epsfig,amsthm,a4,amsmath,url}
\usepackage[latin1]{inputenc}

\newtheorem{Thm}{Theorem}[section]
\newtheorem{Prop}[Thm]{Proposition}
\newtheorem{Lem}[Thm]{Lemma}
\newtheorem{Cor}[Thm]{Corollary}
\theoremstyle{definition}
\newtheorem{Def}[Thm]{Definition}

\newtheorem{Rem}[Thm]{Remark}

\newtheorem{Conj}{Conjecture}

\newcommand{\Z}{\mathbf{Z}}

\newcommand{\N}{\mathbf{N}}
\newcommand{\R}{\mathbf{R}}

\newcommand{\Q}{\mathbf{Q}}

\newcommand{\Isom}{\text{Isom}}

\newcommand{\HH}{\mathcal{H}}
\newcommand{\bpr}{\noindent \textbf{Proof}: }

\newcommand{\eps}{\varepsilon}
\newcommand{\ep}{\hfill$\square$}

\title{Isometric group actions on Hilbert spaces:
growth of cocycles}
\author{Yves de Cornulier, Romain Tessera, Alain Valette}
\date{\today}

\begin{document}

\baselineskip=16pt

\maketitle

\begin{abstract}
We study growth of 1-cocycles of locally compact groups, with values
in unitary representations. Discussing the existence of 1-cocycles
with linear growth, we obtain the following alternative for a class
of amenable groups $G$ containing polycyclic groups and connected
amenable Lie groups: either $G$ has no quasi-isometric embedding
into a Hilbert space, or $G$ admits a proper cocompact action on
some Euclidean space.

On the other hand, noting that almost coboundaries (i.e. 1-cocycles
approximable by bounded 1-cocycles) have sublinear growth, we
discuss the converse, which turns out to hold for amenable groups
with ``controlled" F\o lner sequences; for general amenable groups
we prove the weaker result that 1-cocycles with sufficiently small
growth are almost coboundaries. Besides, we show that there exist,
on a-T-menable groups, proper cocycles with arbitrary small growth.

\medskip

    \hfill\break
\noindent {\sl Mathematics Subject Classification:} Primary 22D10;
Secondary 43A07, 43A35, 20F69. \hfill\break {\sl Key words and
Phrases:} Haagerup property, a-T-menability, amenability, growth of
cocycles, Hilbert distances, geometric group theory, Bernstein
functions.
\end{abstract}


\noindent \textbf{Notation.} Let $G$ be a locally compact group, and
$f,g:G\to\R_+$. We write $f\preceq g$ if there exists $M>0$ and a
compact subset $K\subset G$ such that $f\le Mg$ outside $K$. We
write $f\sim g$ if $f\preceq g\preceq f$. We write $f\prec g$ if,
for every $\eps>0$, there exists a compact subset $K\subset G$ such
that $f\le\eps g$ outside $K$.

\section{Introduction}

The study of affine isometric actions on Hilbert spaces has proven
to be a fundamental tool in geometric group theory. Let $G$ be a
locally compact group, and $\alpha$ an affine isometric action on an
affine Hilbert space $\HH$ (real or complex). The function
$b:G\to\HH$ defined by $b(g)=\alpha(g)(0)$ is called a 1-cocycle
(see Section \ref{Sec:growth} for details), and we call the function
$g\mapsto\|b(g)\|$ a \textit{Hilbert length function} on~$G$. We
focus on the \textit{asymptotic} behaviour of Hilbert length
functions on a given group $G$. A general question is the following:
how is it related to the structure of~$G$?

For instance, if $G$ is $\sigma$-compact, $G$ has the celebrated
Kazhdan's Property (T) if and only if every Hilbert length function
is bounded (see \cite{HV}). This is known to have strong
group-theoretic consequences on~$G$: for instance, this implies that
$G$ is compactly generated and has compact abelianization (see
\cite[Chap. 2]{BHV} for a direct proof).

In this paper, we rather deal with groups which are far from having
Kazhdan's Property (T): a locally compact group $G$ is called
a-T-menable if it has a proper Hilbert length function. The class of
a-T-menable locally compact groups contains (see \cite{CCJJV})
amenable groups, Coxeter groups, isometry groups of locally finite
trees, isometry groups of real and complex hyperbolic spaces and all
their closed subgroups, such as free and surface groups. We show in
§\ref{Sec:slow_growth} that, for a-T-menable locally compact groups
(for instance, $\Z$), there exist proper Hilbert length functions of
arbitrary slow growth.

The study of Hilbert length functions with non-slow growth is more
delicate. An easy but useful observation is that, for a given
compactly generated, locally compact group, any Hilbert length
function $L$ is linearly bounded, i.e. $L(g)\preceq |g|_S$, where
$|\cdot|_S$ denotes the word length with respect to some compact
generating subset.

We discuss, in Section \ref{Sec:sublin_growth}, Hilbert length
functions with sublinear growth. These include those Hilbert length
functions whose corresponding 1-cocycle (see Section
\ref{Sec:growth}) is an almost coboundary, i.e. can be approximated,
uniformly on compact subsets, by bounded 1-cocycles. We discuss the
converse.

Denote by $(\mathcal{L})$ the class of groups including:
\begin{itemize}
\item polycyclic groups and connected amenable Lie groups,

\item semidirect products
    $\Z[\frac{1}{mn}]\rtimes_{\frac{m}{n}}\Z$, with
    $m,n$ co-prime integers with $|mn|\ge 2$ (if $n=1$ this is the
    Baumslag-Solitar group $BS(1,m)$); semidirect products
    $\left(\R\oplus\bigoplus_{p\in S}\Q_p\right)\rtimes_{\frac{m}{n}}\Z$ or
    $\left(\bigoplus_{p\in S}\Q_p\right)\rtimes_{\frac{m}{n}}\Z$, with
    $m,n$ co-prime integers, and $S$ a finite set of prime numbers
    dividing $mn$.

\item wreath products $F\wr\Z$ for $F$ a finite group.
\end{itemize}

\begin{Thm}[\textnormal{see Corollary \ref{controlfol},
Propositions \ref{prop:amenable_slow_implies_alm_cob} and
\ref{prop:cocompact_SOn1}}]~\begin{itemize} \item [(1)] If $G$ is a
compactly generated, locally compact amenable group, then every
1-cocycle with sufficiently slow growth is an almost boundary.

\item [(2)] For groups in the class $(\mathcal{L})$, every
sublinear 1-cocycle is an almost coboundary.

\item [(3)] If $\Gamma$ is a finitely generated, discrete, non-amenable subgroup in
$\textnormal{SO}(n,1)$ or $\textnormal{SU}(n,1)$ for some $n\ge 2$,
then $\Gamma$ admits a 1-cocycle with sublinear growth (actually
$\preceq |g|^{1/2}$) which is not an almost coboundary.
\end{itemize}
\end{Thm}

In §\ref{subs:large_sublinear}, we show that there exist, on $\R^n$
or $\Z^n$, Hilbert length functions with arbitrary large sublinear
growth, showing that, in a certain sense, there is no gap between
Hilbert length functions of linear and of sublinear growth.

\medskip

In Section \ref{Sec:lin_growth}, we discuss the existence of a
Hilbert length function on $G$ with linear growth. Such a function
exists when $G=\Z^n$. We conjecture that the converse is
essentially true.

\begin{Conj}\label{conj1} Let $G$ be a locally compact, compactly generated
group having a Hilbert length function with linear growth. Then $G$
has a proper, cocompact action on a Euclidean space. In particular,
if $G$ is discrete, then it must be virtually abelian.
\label{Conj_INTRO_cocycle_linear}
\end{Conj}

Our first result towards Conjecture \ref{Conj_INTRO_cocycle_linear}
is a generalization of a result by Guentner and Kaminker
\cite[§5]{GK} to the non-discrete case.

\begin{Thm}[\textnormal{see Theorem \ref{thm:GK_lc}}]
Let $G$ be a locally compact, compactly generated group. If $G$
admits a Hilbert length function with growth $\succ |g|^{1/2}$ (in
particular, if it admits a Hilbert length function with linear
growth), then $G$ is amenable.
\end{Thm}

We actually provide a new, simpler proof, while it is not clear how
to generalize the proof in \cite{GK} to the non-discrete case
\footnote{The truncation argument in \cite{GK} does not seem to
extend to non-discrete groups.}.

To prove that locally compact groups in the class $(\mathcal{L})$
satisfy Conjecture \ref{Conj_INTRO_cocycle_linear}, we use Shalom's
Property $H_{FD}$: a locally compact group has Property $H_{FD}$ if
any unitary representation with nontrivial reduced cohomology has a
finite-dimensional nonzero subrepresentation. All groups in the
class $(\mathcal{L})$ are known to satisfy Property $H_{FD}$. We
prove

\begin{Thm}[\textnormal{see Theorem \ref{FHDcoclineaires}}]Locally compact,
compactly generated groups with Property $H_{FD}$ satisfy Conjecture
\ref{Conj_INTRO_cocycle_linear}.
\end{Thm}

We next consider uniform embeddings into Hilbert spaces. There is a
nice trick, for which we are indebted to Gromov, allowing to
construct, if the group is amenable, a 1-cocycle with the same
growth behaviour as the initial embedding. See
Proposition~\ref{prop:rem_gromov} for a precise statement. Thus we
get:

\begin{Thm}If $G$ is any locally compact, compactly generated, amenable group with Property $H_{FD}$
(e.g. in the class ($\mathcal{L})$), then
\begin{itemize}
    \item
either $G$ does \textit{not} admit any quasi-isometric embedding
into a Hilbert space,

\item or $G$ acts properly cocompactly on some Euclidean space
(i.e. a finite-dimensional real Hilbert space).

\end{itemize}\label{thm:intro_noQIemb}\end{Thm}

The reader interested in the proof of Theorem
\ref{thm:intro_noQIemb} can skip Section \ref{Sec:sublin_growth},
except the elementary Proposition \ref{coboundaries}. Let us observe
that the proof of Theorem \ref{thm:intro_noQIemb} does {\it not}
appeal to asymptotic cones. It contains, as a particular case, the
fact that a simply connected nilpotent non-abelian Lie group has no
quasi-isometric embedding into a Hilbert space, a result due to S.
Pauls \cite{Pauls}. Moreover, Theorem \ref{thm:intro_noQIemb}
provides new proofs of two known results (see §\ref{subs:unif_emb}
for proofs):

\begin{Cor} [Quasi-isometric rigidity of $\Z^n$]
If a finitely generated group is quasi-isometric to $\Z^n$, then it
has a finite index subgroup isomorphic to
$\Z^n$.\label{cor:rigidite_virt_ab}
\end{Cor}

This latter result has recently been proved by Shalom
\cite{Shalom04}, who first establishes the invariance of Property
$H_{FD}$ by quasi-isometries. We also make use of this crucial fact,
although the use of Proposition \ref{prop:rem_gromov} allows us to
conclude in a different way.

\begin{Cor}[Bourgain \cite{Bourgain}]
For $r\geq 3$, the regular tree of degree $r$ does not embed
quasi-isometrically into a Hilbert
space.\label{cor:tree_not_QI_embed_Hilb}
\end{Cor}

A locally compact, compactly generated group is either non-amenable
or non-unimodular if and only if it is quasi-isometric to a graph
with positive Cheeger constant (see~\cite[Corollary~7.12]{tess}). In
Corollary 1.6 of \cite{BenSc}, Benjamini and Schramm use Bourgain's
result above, to prove that a graph with positive Cheeger constant
cannot be quasi-isometrically embedded into a Hilbert space. As a
consequence: {\it if a compactly generated, locally compact group
$G$ admits a quasi-isometric embedding into Hilbert space, then $G$
is amenable and unimodular}.

By Proposition \ref{prop:rem_gromov}, the existence of a
quasi-isometric embedding into a Hilbert space implies, for an
amenable group, the existence of a Hilbert length function with
linear growth. So Conjecture~\ref{conj1} is equivalent to the
following statement, apparently more general:

\noindent {\bf Conjecture 1'} A locally compact, compactly generated
group $G$ admitting a quasi-isometric embedding into a Hilbert space
has a proper, cocompact action on a Euclidean space. In particular,
if $G$ is discrete, then it must be virtually abelian.

\medskip

We conclude this introduction with a remark about compression. The
following definition is due to E. Guentner and J. Kaminker
\cite{GK}. Let $G$ be a compactly generated group, endowed with its
word length $|.|_{S}$.

\begin{Def} The {\em equivariant Hilbert space compression} of $G$ is
defined as: $$B(G)=\sup\{\alpha\ge 0,\;\exists\textnormal{ unitary
representation }\pi\;\;\exists b\in Z^1(G,\pi),\; \|b(g)\|\succeq
|g|_S^\alpha\}.$$
\end{Def}

It is clear that $0\le B(G)\le 1$, and if $G$ admits a linear
1-cocycle, then $B(G)=1$. The converse in not true: it is shown in
\cite{Tes} that $B(G)=1$ for groups of the class $(\mathcal{L})$
whereas we have shown above (Theorem \ref{thm:intro_noQIemb}) that
these groups do not admit linear cocycles unless they act properly
on a Euclidean space.

Another immediate observation is that if $B(G)>0$, then $G$ is
a-T-menable. We know nothing about the converse: actually we know no
example of an a-T-menable group with $B<1/2$; at the other extreme,
we do not know if solvable groups always satisfy $B>0$.

It follows from Proposition \ref{prop:rem_gromov} that, for
\textit{amenable groups} $G$, {\it the number $B(G)$ is a
quasi-isometry invariant}. This probably does not hold for
non-amenable groups, but we do not know any counterexample. More
precisely

\begin{itemize}
\item It is not known if being a-T-menable is a quasi-isometry invariant.
\item It is not known if there exists a non-amenable, a-T-menable,
compactly generated, locally compact group $G$ with $B(G)\neq 1/2$.
\end{itemize}

Finally, let us mention that it can be interesting to study the
growth of 1-cocycles when we restrict to certain special classes of
unitary representations. In the case of the regular representation,
some results can be found in \cite{Tes}, related to isoperimetric
properties of the group.

\bigskip

\noindent \textbf{Acknowledgments.} We are indebted to Misha Gromov
for a decisive remark. We also thank Emmanuel Breuillard, Pierre de
la Harpe, and Urs Lang for useful remarks and corrections.

\section{Preliminaries}\label{Sec:growth}

\subsection{Growth of 1-cocycles}

Let $G$ be a locally compact group, and $\pi$ a unitary or
orthogonal representation (always assumed continuous) on a Hilbert
space $\mathcal{H}=\mathcal{H}_{\pi}$. The space $Z^1(G,\pi)$  is
defined as the set of continuous functions $b:G\to\mathcal{H}$
satisfying, for all $g,h\in G$, the 1-cocycle condition
$b(gh)=\pi(g)b(h)+b(g)$. Observe that, given a continuous function
$b:G\to\mathcal{H}$, the condition $b\in Z^1(G,\pi)$ is equivalent
to saying that $G$ acts by affine transformations on $\mathcal{H}$
by $\alpha(g)v=\pi(g)v+b(g)$. The space $Z^1(G,\pi)$ is endowed with
the topology of uniform convergence on compact subsets.

The subspace of coboundaries $B^1(G,\pi)$ is the subspace (not
necessarily closed) of $Z^1(G,\pi)$ consisting of functions of the
form $g\mapsto v-\pi(g)v$ for some $v\in\mathcal{H}$. It is
well-known \cite[§4.a]{HV} that $b\in B^1(G,\pi)$ if and only if $b$
is bounded on $G$.

The subspace of almost coboundaries $\overline{B^1(G,\pi)}$ is the
closure of $B^1(G,\pi)$. A 1-cocycle $b$ is an almost coboundary if
and only if the corresponding affine action \textit{almost has fixed
points}, i.e. for every compact subset $K\subset G$ and $\eps>0$,
there exists $v$ such that $\sup_{g\in K}\|\alpha(g)v-v\|\le\eps$
(see \cite[§3.1]{BHV}). When $G$ is generated by a symmetric compact
subset $S$, it suffices to check this condition for $K=S$, and a
sequence of \textit{almost fixed points} is defined as a sequence
$(v_n)$ such that $\sup_{g\in S}\|\alpha(g)v_n-v_n\|\to 0$.

The first cohomology space of $\pi$ is defined as the quotient space
$H^1(G,\pi)=Z^1(G,\pi)/B^1(G,\pi)$, and the first reduced cohomology
space of $\pi$ is defined as
$\overline{H^1}(G,\pi)=Z^1(G,\pi)/\overline{B^1(G,\pi)}$.

\medskip

Now suppose that $G$ is a locally compact, compactly generated
group. For $g\in G$, denote by $|g|_{S}$ the word length of $g$ with
respect to an open, relatively compact generating set $S\subset G$.

Let $b\in Z^1(G,\pi)$ be a 1-cocycle with respect to a unitary
representation $\pi$ of $G$. We study the growth of $\|b(g)\|$ as a
function of~$g$.

\begin{Def}
The \textit{compression} of the 1-cocycle $b$ is the function
$$\rho:\R_{+}\rightarrow \R_{+}\cup\{\infty\}:\;x\mapsto
\rho(x)=\inf\{\|b(g)\|:\;g\in G,\,|g|_S\ge x\}.$$
\end{Def}

\begin{Rem}
A related notion is the \textit{distortion function}, defined in
\cite{Farb} in the context of an embedding between finitely
generated groups. The distortion function of the 1-cocycle $b$ is
defined as the function $\R_+\to\R_+\cup\{\infty\}$ by
$f(x)=\sup\{|g|_{S}\,:\;\|b(g)\|\le x\}$. The reader can check that,
except in trivial cases\footnote{Trivial cases are: when $G$ is
compact, so that $\rho$ is eventually equal to $\infty$ and $f$ is
eventually equal to a finite constant, and when $b$ is not proper,
so that $\rho$ is bounded, and $f$ is eventually equal
to~$\infty$.}, the compression $\rho$ and the distortion $f$ are
essentially reciprocal to each other.
\end{Rem}

Recall that a \textit{length function} on a group $\Gamma$ is a
function $L:\Gamma\to\R_+$ satisfying $L(1)=0$ and, for all
$g,h\in\Gamma$, $L(g^{-1})=L(g)$ and $L(gh)\le L(g)+L(h)$, so that
$d(g,h)=L(g^{-1}h)$ is a left-invariant pseudo-distance ({\it
``\'ecart''}) on $\Gamma$.

It is immediate from the 1-cocycle relation that the function
$g\mapsto \|b(g)\|$ is a length function on the group $G$. In
particular, if $G$ is locally compact, compactly generated, then it
is dominated by the word length. We thus obtain the following
obvious bound:

\begin{Prop}\label{cocycles}
For $b\in Z^1(G,\pi)$, we have $\|b(g)\|\preceq |g|_{S}$.\ep
\end{Prop}

\medskip

Define $$\text{lin}(G,\pi)=\{b\in Z^1(G,\pi),\;\|b(g)\|\succeq
|g|_S\}$$
$$\text{sublin}(G,\pi)=\{b\in Z^1(G,\pi),
\|b(g)\|\prec |g|_S\},$$ namely, the set of cocycles with linear
(resp. sublinear) growth. Here are immediate observations:

\begin{itemize}

\item $\text{sublin}(G,\pi)$ is a linear subspace of $Z^1(G,\pi)$.

\item $B^1(G,\pi)\;\subset\; \text{sublin}(G,\pi)\;\subset\;
Z^1(G,\pi)\smallsetminus \text{lin}(G,\pi)$.

\item If $G=\Z$ or $\R$, then it is easy to check that
$Z^1(G,\pi)=\text{lin}(G,\pi)\cup\text{sublin}(G,\pi)$ (this follows
either from Corollary \ref{controlfol} below, or from a direct
computation
\footnote{Here is the argument for $\Z$. Fix $b\in Z^1(\Z,\pi)$.
Write $U=\pi(1)$ and $\pi=\pi_0\oplus\pi_1$, where
$\pi_0=\text{Ker}(U-1)$ denotes the invariant vectors, and decompose
$b$ as $b_0+b_1$. Clearly, $b_0$ is either zero or has linear
growth. On the other hand, $b_1$ has sublinear growth: indeed, as
$b_1(1)$ is orthogonal to $\text{Ker}(U-1)$, it belongs to
$\overline{\text{Im}(U-1)}$, so that $b_1\in
\overline{B^1(\Z,\pi_1)}$. Hence $b_1\in\text{sublin}(\Z,\pi)$ by
Corollary \ref{B1bar} below.}).
On the other hand, this does not generalize to arbitrary~$G$.
Indeed, take any nontrivial action of $\Z^2$ by translations on
$\R$: then the associated cocycle is neither linear nor sublinear.

\end{itemize}

\subsection{Conditionally negative definite functions and Bernstein
functions}\label{subsbernstein}

A conditionally negative definite function on a group $G$ is a
function $\psi:G\to\R_+$ such that $\psi^{1/2}$ is a Hilbert length
function. Equivalently \cite[5.b]{HV}, $\psi(1)=0$,
$\psi(g)=\psi(g^{-1})$ for all $g$, and, for all
$\lambda_1,\dots,\lambda_n\in\R$ such that $\sum_{i=1}^n\lambda_i=0$
and for all $g_1,\dots,g_n\in G$, we have
$\sum_{i,j=1}^n\lambda_i\lambda_j\psi(g_i^{-1}g_j)\le 0$. Continuous
conditionally negative definite functions on a locally compact group
$G$ form a convex cone, closed under the topology of uniform
convergence on compact subsets.

\medskip

A continuous function $F:\R_{+}\rightarrow\R_{+}$ is a {\it
Bernstein function} if there exists a positive measure $\mu$ on
Borel subsets of $\R^{*}_+$ such that $\mu([\eps,\infty[)<\infty$
for all $\eps>0$, $\int_0^1x\,d\mu(x)<\infty$, and such that, for
some $a\ge 0$,
$$\forall t>0,\quad F(t)\,=\,at\,+\,\int_{0}^{+\infty} (1-e^{-tx})\,d\mu(x).$$

Note that such a function is real analytic on $\R_+^*$. We note for
reference the following well-known result due to of Bochner and
Schoenberg \cite[Theorem 8]{Sch}:

\begin{Lem}Let $\psi$ be a conditionally negative definite function on $G$,
and let $F$ be a Bernstein function. Then $F\circ\psi$ is
conditionally negative definite on
$G$.\ep\label{lem:Bochner_Schoenberg}
\end{Lem}

Examples of Bernstein functions are $x\mapsto x^a$ for $0< a\le 1$,
and $x\mapsto\log(x+1)$. For more on Bernstein functions, see for
instance \cite{BF}.

\section{Cocycles with sublinear growth}\label{Sec:sublin_growth}

\subsection{Almost coboundaries are sublinear}

\begin{Prop}\label{coboundaries}
Let $G$ be a locally compact, compactly generated group. In
$Z^1(G,\pi)$, endowed with topology of uniform convergence on
compact subsets,

1) $\textnormal{sublin}(G,\pi)$ is a closed subspace;

2) $\textnormal{lin}(G,\pi)$ is an open subset.
\end{Prop}

{\bf Proof:} Fix a symmetric open, relatively compact generating
subset $S\subset G$. Let $b$ be the limit of a net $(b_i)_{i\in I}$
in $Z^1(G,\pi)$. Write $b'_i=b-b_i$, and fix $\eps>0$. For $i$ large
enough (say, $i\ge i_0$), $\sup_{s\in S}\|b'_{i}(s)\|\le\eps/2$.
Since $g\mapsto\|b'_{i}(g)\|$ is a length function, this implies
that for every $g\in G$ and $i\geq i_{0}$,
$\|b'_i(g)\|\le\eps|g|_S/2$, i.e. $\|b'_i(g)\|/|g|_{S}\le\eps/2$.

\begin{enumerate}
    \item  [1)] Suppose that all $b_i$'s belong to $\text{sublin}(G,\pi)$. Fix
$i\geq i_{0}$. Then $\|b_{i}(g)\|/|g|_{S}\leq \eps/2$ for $g$ large
enough (say, $g\notin K$ compact). So $\|b(g)\|/|g|_{S}\le\eps$ for
$g\notin K$. This shows that $b\in \text{sublin}(\Gamma,\pi)$, so we
are done.

    \item  [2)] Suppose that $b\in \text{lin}(G,\pi)$. Then, if $\eps$ has been
chosen sufficiently small, $\|b(g)\|/\|g\|\ge\eps$ for large $g$
(say, $g\notin K$ compact). Hence,
$\|b_i(g)\|/|g|_S\ge(\|b(g)\|-\|b'_i(g)\|)/|g|_S\ge \eps/2$ for
$i\ge i_0$, $g\notin K$, showing that $b_i\in \text{lin}(G,\pi)$ for
$i\ge i_0$. \hfill$\square$
\end{enumerate}

\begin{Rem}
In the previous result, it is essential that we fix the unitary
representation $\pi$. Indeed, it is easy to show that, on every
group $G$, every Hilbert length function (e.g. of linear growth) can
be approximated, uniformly on compact subsets, by bounded Hilbert
length functions.\end{Rem}

\begin{Cor}\label{B1bar}
If $b\in \overline{B^1}(G,\pi)$, then $\|b(g)\|\prec
|g|_{S}$.\label{cor:Z1bar_sslineaire}
\end{Cor}

{\bf Proof:} $B^1(G,\pi)\subset \text{sublin}(G,\pi)$, so that
$$\overline{B^1(G,\pi)}\subset
\overline{\text{sublin}(G,\pi)}=\text{sublin}(G,\pi)$$ by
Proposition \ref{coboundaries}. \hfill$\square$

\subsection{Groups with controlled F\o lner sequences}

In this section, we prove that the converse to Corollary
\ref{cor:Z1bar_sslineaire} is true for unimodular groups in the
class $(\mathcal{L})$: that is, a cocycle has sublinear growth if
and only if it is an almost coboundary.

\medskip

Let $G$ be a compactly generated, locally compact group with Haar
measure $\mu$, and let $S$ be a compact generating subset. Let
$(F_n)$ be a sequence of measurable, bounded subsets of nonzero
measure. Set
$$\eps_n=\frac{\sup_{s\in S}\mu(sF_n\triangle F_n)}{\mu(F_n)}.$$

Consider an isometric affine action $\alpha$ of $G$ on a Hilbert
space, and let $b$ be the corresponding 1-cocycle. Set
$$v_n=\frac{1}{\mu(F_n)}\int_{F_n}b(g)d\mu(g).$$
This is well-defined.

\begin{Lem}
Suppose that $\sup_{g\in F_n}\|b(g)\|\prec 1/\eps_n$. Then $(v_n)$
is a sequence of almost fixed points for the affine action $\alpha$
associated with $b$.\label{lem:moyenne_Folner_almost_inv}
\end{Lem}
\bpr For $s\in S$, we have
$$\alpha(s)v_n-v_n=\frac{1}{\mu(F_n)}\int_{F_n}(b(sg)-b(g))d\mu(g).$$
Thus
$$\|\alpha(s)v_n-v_n\|\;\le\;\frac{2}{\mu(F_n)}\int_{sF_n\triangle F_n}
\|b(g)\|d\mu(g)\;\le\; 2\eps_n\sup_{g\in F_n}\|b(g)\|.\;\;\Box$$

\medskip

Recall (see \cite[Appendix G]{BHV}) that ``$G$ is amenable'' exactly
means that we can choose $(F_n)$ so that $\eps_n\to 0$, and $(F_n)$
is then called a F\o lner sequence. In this case, we obtain, as a
consequence of Lemma \ref{lem:moyenne_Folner_almost_inv}, that a
1-cocycle of sufficiently slow growth (depending on the behaviour of
the F\o lner sequence, i.e. on the asymptotic behaviour of $\eps_n$
and the diameter of $(F_n)$) must be an almost coboundary. We record
this as:

\begin{Prop}Let $G$ be a compactly generated, locally compact
amenable group. Then there exists a function
$u:G\to\R_+\cup\{\infty\}$ satisfying
\begin{itemize}
\item $\lim_{g\to \infty}u(g)=\infty$;

\item for every 1-cocycle $b$ of $G$, $\|b(g)\|\prec u(g)$ implies
that $b$ is an almost coboundary.
\end{itemize}
In other words, if two 1-cocycles are sufficiently close, then they
coincide in reduced
1-cohomology.\ep\label{prop:amenable_slow_implies_alm_cob}
\end{Prop}

Explicitly, the function $u$ can be defined as follows
$$u(g)=\frac{1}{\max\{\eps_n|\;n\in\N\text{ s.t. }g\in F_n\}},$$
where we set $u(g)=\infty$ if $g\notin\bigcup F_n$.

To obtain stronger statements we introduce a more restrictive notion
of F\o lner sets.

\begin{Def} We say that the F\o lner sequence $(F_n)$ of the amenable,
compactly generated locally compact group $G$ is {\it controlled} if
there exists a constant $c>0$ such that, for all~$n$,
$$F_n\subset B(1,c/\eps_n).$$
\end{Def}

In \cite{Tes}, it is proved that a unimodular group in the class
$(\mathcal{L})$ admits a controlled F\o lner sequence.

\begin{Cor}\label{controlfol}
Let $G$ be a compactly generated, locally compact amenable group
admitting a controlled F\o lner sequence $(F_n)$, and keep the
notation as above. Then the following statements are equivalent:
 \\(1) $b\in\overline{B^1}(\Gamma,\pi)$
 \\(2) $b\in\textnormal{sublin}(\Gamma,\pi)$
 \\(3)The sequence $(v_n)$ is a sequence of almost fixed points for
$\alpha$.\ep
\end{Cor}

{\bf Proof:} (3)$\Rightarrow$(1) is immediate , while
(1)$\Rightarrow$(2) follows from Corollary
\ref{cor:Z1bar_sslineaire}. The remaining implication is
(2)$\Rightarrow$(3): suppose that $b$ is sublinear. Write
$$f(r)=\sup_{|g|\le r}\|b(g)\|,$$ where $f(r)\prec r$. Then
$$\sup_{g\in F_n}\|b(g)\|\;\le\; \sup_{|g|\le
c/\eps_n}\|b(g)\|\;=\;f(c/\eps_n)\;\prec\; 1/\eps_n,$$ so that we
can apply Lemma \ref{lem:moyenne_Folner_almost_inv} to obtain that
$(v_n)$ is a sequence of almost invariant vectors.\ep

We use this to prove a conjecture of Shalom \cite[Section
6]{Shalom04}. Recall that a representation of a group $\Gamma$ is
said to be {\em finite} if it factors through a finite group.

\begin{Prop}
Let $\pi$ be a unitary representation of a finitely generated,
virtually nilpotent group $\Gamma$ and let $S$ be a finite
generating subset of $\Gamma$. Suppose that $\pi$ has no finite
subrepresentation\footnote{In the conjecture of \cite{Shalom04} the
assumptions are slightly stronger: $S$ is assumed symmetric, and
$\pi$ is supposed to have no finite-dimensional subrepresentation.}.
For every cocycle $b\in Z^1(\Gamma,\pi)$, define:
$$v_n=\frac{1}{|S^n|}\sum_{g\in S^n}b(g).$$
Then there exists a subsequence $(v_{n_i})$ which is a sequence of
almost fixed points for the affine action $\alpha$ associated with
$b$:
$$\|\alpha(s)v_{n_i}-v_{n_i}\|\rightarrow 0,\;\forall s\in
S.$$\label{thm:conj_shalom}
\end{Prop}

{\bf Proof:} First recall \cite{Wolf} that there exists $d>0$ such
that $|S^n|\preceq n^d$ for all~$n$. By an elementary
argument\footnote{If $(|S^{n+1}|-|S^{n}|)/|S^{n}|\ge d/n$ for all
(large) $n$, then $|S^n|\succeq n^d$.}, there exists an infinite
sequence $(n_i)$ such that:
\begin{equation}\label{eq:sphere}
\frac{|S^{n_i+1}\setminus
S^{n_i}|}{|S^{n_i}|}\preceq\frac{1}{|n_i|}.
\end{equation}
It follows that the family $(S^{n_i})_{i}$ is a controlled F\o lner
sequence of $\Gamma$.

Since $\Gamma$ is virtually nilpotent, by Corollary 5.1.3 and Lemma
4.2.2 in \cite{Shalom04}, it has property $H_F$, i.e. every
representation with non-zero first reduced cohomology has a finite
subrepresentation. Here, by our assumption:
$\overline{H^1}(\Gamma,\pi)=0$. So the conclusion follows from
Corollary \ref{controlfol}. \hfill$\square$
\medskip

{\bf Remark:} Shalom proved in the final section of \cite{Shalom04}
that, if the result of Proposition \ref{thm:conj_shalom} was proved
under the bare assumption that $\Gamma$ has polynomial growth, this
would give rise to a new, simpler\footnote{The proof would be
simpler in that it would not appeal to the solution of Hilbert's 5th
problem about the structure of locally compact groups.} proof of
Gromov's celebrated theorem \cite{Gromov}: a finitely generated
group of polynomial growth is virtually nilpotent.

\hyphenation{cohomo-logy}

\subsection{A sublinear cocycle with nontrivial reduced 1-cohomo\-logy}

It turns out that the converse of Corollary
\ref{cor:Z1bar_sslineaire} is not true in general, for finitely
generated groups.

\begin{Prop} Let $\Gamma$ be a discrete, finitely generated, nonamenable subgroup either in
$G=\textnormal{SO}(n,1)$ $(n\ge 2)$ or $G=\textnormal{SU}(m,1)$
$(m\geq 1)$. There exists a unitary representation $\sigma$ of
$\Gamma$, and $b\in
Z^1(\Gamma,\sigma)-\overline{B^1}(\Gamma,\sigma)$, such that
$$\|b(g)\|\preceq |g|_{S}^{1/2}.$$
If moreover $\Gamma$ is a cocompact lattice and either $n\geq 3$ or
$m\geq 2$, then $\|b(g)\|\sim |g|_S^{1/2}$ and the representation
$\sigma$ may be taken to be irreducible.\label{prop:cocompact_SOn1}
\end{Prop}

{\bf Proof:} A result of Delorme \cite[Lemme~V.5]{Del} says that
there exists a unitary irreducible representation $\pi$ with
$\overline{H^1}(G,\pi)\neq 0$: so we choose $b\in
Z^1(G,\pi)-\overline{B^1}(G,\pi)$. Let $K$ be a maximal compact
subgroup of $G$; replacing $b$ by a cohomologous 1-cocycle, we may
assume that $b|_{K}\equiv 0$. Then $b:G\rightarrow{\mathcal
H}_{\pi}$ factors through a map $F:G/K\rightarrow{\mathcal
H}_{\pi}$, which is equivariant with respect to the corresponding
affine action on ${\mathcal H}_{\pi}$. By an unpublished result of
Shalom (for a proof, see Corollary 3.3.10 in \cite{BHV}), the map
$F$ is harmonic. We may now appeal to Gromov's results
(\cite{GroRandom}, section 3.7.D'; see also Proposition 3.3.21 in
\cite{BHV}) on the growth of harmonic, equivariant maps from a rank
1, Riemannian, symmetric space to a Hilbert space. If $d(x,x_{0})$
denotes the Riemannian distance between $x$ and the point $x_{0}$
with stabilizer $K$ in $G/K$, then for some constant $C>0$,
$$\|F(x)\|^2\,=\,C\,d(x,x_{0}) + o(d(x,x_{0})).$$
Set $\sigma=\pi|_{\Gamma}$; then
$$\|b|_{\Gamma}(g)\|^2\,\sim\,d(gx_0,x_0) \preceq |g|_{S}.$$

Let us prove that $b|_{\Gamma}$ is not an almost coboundary. First,
we appeal to a result of Shalom \cite{ShaAnn}, stating that the
restriction map $H^1(G,\pi)\rightarrow H^1(\Gamma,\sigma)$ is
injective, so that $b|_{\Gamma}$ is unbounded. Now $\sigma$, as the
restriction to $\Gamma$ of a non-trivial unitary irreducible
representation of $G$, does not weakly contain the trivial
representation of $\Gamma$ (this follows from a result of Cowling
\cite{Cow}: for some $p>1$, all coefficients of $\pi$ are in
$L^p(G)$. Therefore, a suitable tensor power of $\pi$ is a
sub-representation of $\lambda_G \oplus \lambda_G \oplus ...$, where
$\lambda_G$ is the regular representation; restricting, a suitable
tensor power of $\sigma$ is a sub-representation of
$\lambda_{\Gamma}\oplus \lambda_{\Gamma}\oplus...$. Non-amenability
of $\Gamma$ allows one to conclude). By Guichardet's well-known
criterion (see \cite{Guich}, Cor. 2.3 in Chap. III), this implies
that $B^1(\Gamma,\sigma)=\overline{B^1(\Gamma,\sigma)}$, in
particular every almost coboundary is bounded. Since $b|_{\Gamma}$
is unbounded, this gives the desired result.

Finally, if $\Gamma$ is a cocompact lattice in $G$, then $\Gamma$ is
quasi-isometric to $G/K$, so we get
$\|b|_{\Gamma}(g)\|^2\,\sim\,|g|_{S}.$ Moreover, if $n\geq 3$ or
$m\geq 2$, then $\pi$ is not in the discrete series of $G$ (see
\cite[Remarque V.8]{Del}), so that $\sigma=\pi|_{\Gamma}$ is
irreducible, by a result of Cowling and Steger \cite[Prop.
2.5]{CowSte}. \hfill$\square$

\subsection{Cocycles with slow growth}\label{Sec:slow_growth}

We prove here that, on an a-T-menable group (e.g. $\Z$), there exist
cocycles with arbitrarily slow growth.

\begin{Prop}\label{slowgrowth}
    Assume that $G$ is locally compact, a-T-menable. For every proper
    function $f:G\rightarrow [1,\infty[$, there
    exists a continuous conditionally negative definite, proper function $\psi$
    on $G$ such that $\psi\leq f$.
\end{Prop}

This is obtained as a consequence of the following lemma (see
\S\ref{subsbernstein} for the definition of Bernstein functions).

\begin{Lem}\label{Bernstein} Let $u$ be a proper function on
$\R_{+}$, with $u(t)\geq 1$ for $t\in\R_+$. There exists a proper
Bernstein function $F$ such that $F(t)\leq u(t)$ for $t\in\R_+$.
\end{Lem}

{\bf Proof:} We are going to define inductively a sequence
$(x_{n})_{n\geq 1}$ of positive real numbers such that
$0<x_{n}<2^{-n}$, and define
$$F(t)=\sum_{n=1}^{\infty}(1-e^{-tx_{n}}).$$
Since $1-e^{-tx_{n}}\leq tx_{n}$, the series defining $F$ will
converge uniformly on compact subsets of $\R_{+}$, so $F$ will be a
Bernstein function (in fact associated with
$\mu=\sum_{n=1}^{\infty}\delta_{x_{n}}$ and $a=0$). Let
$F_{m}(t)=\sum_{n=1}^{m}(1-e^{-tx_{n}})$ be the $m$-th partial sum.
For fixed $m$, we will have $F\geq F_{m}$, hence
$$\liminf_{t\rightarrow\infty}F(t)\,\geq\,\lim_{t\rightarrow\infty}F_{m}(t)
\,=\,m;$$ since this holds for every $m$, we have
$\lim_{t\rightarrow\infty}F(t)=\infty$, i.e. $F$ is proper.

It remains to manage to construct the $x_{n}$'s so that $F\leq u$ on
$\R_+$. We will construct $x_{n}$ inductively so that
$u>F_{n}+2^{-n}$ on $\R_+$. Setting $F_{0}\equiv 0$, the
construction will also apply to $n=1$. So assume
$0<x_{n-1}<2^{-n+1}$ has been constructed so that
$u>F_{n-1}+2^{-n+1}$ on $\R_+$. Since $u$ is proper and $F_{n-1}$ is
bounded, we find $k_{n}>0$ large enough so that $u(t)>F_{n-1}(t)+2$
for $t>k_{n},\,t\in\R_+$. By taking $x_{n}>0$ very small (with
$x_{n}<2^{-n}$ anyway), we may arrange to have
$1-e^{-tx_{n}}<2^{-n}$ for $t<k_{n}$. Then, for
$t<k_{n},\,t\in\R_+$:
$$u(t)-F_{n}(t)=u(t)-F_{n-1}(t)-(1-e^{-tx_{n}})>
2^{-n+1}-2^{-n}=2^{-n};$$ while for $t\geq k_{n},\,t\in\R_+$:
$$u(t)-F_{n}(t)=u(t)-F_{n-1}(t)-(1-e^{-tx_{n}})>
2-1=1>2^{-n}.$$ This concludes the induction step. \hfill$\square$
\medskip

{\bf Proof of Proposition \ref{slowgrowth}:} If $G$ is compact, we
can take $\psi=0$; thus suppose $G$ noncompact. As $G$ is
a-T-menable, we may choose a proper conditionally negative definite
function $\psi_{0}$ on $G$.

Define a proper function $u\ge 1$ on $\R_+$ by
$$u(t)=\inf\{f(g):\;g\in\psi_0^{-1}(\,[t,\infty[\,)\}.$$
By lemma \ref{Bernstein}, we can find a proper Bernstein function
$F$ such that $F\leq u$ on $\R_+$. Then, by construction,
$F(\psi_{0}(g))\leq f(g)$, and by Lemma
\ref{lem:Bochner_Schoenberg}, $F\circ\psi_{0}$ is conditionally
negative definite. \hfill$\square$

\subsection{Cocycles with arbitrary large sublinear growth}\label{subs:large_sublinear}

As we observed earlier, a cocycle on $\Z^n$ (or $\R^n$) has either
linear or sublinear growth. This raises the question whether there
is a gap between the two. We show here that it is not the case.

\begin{Lem}\label{largesublin}
 Let $w:\R_+\to\R_+$ be any function with sublinear growth. Then there
exists a sublinear Bernstein function $F$ such that $F(x)\ge w(x)$
for $x$ large enough.
 \end{Lem}

\bpr The function $x\mapsto w(x)/x$ tends to zero. It is easy to
construct a decreasing function $x\mapsto u(x)$ of class $C^1$, such
that $u(x)\ge w(x)/x$ for $x$ large enough, and such that $u(x)\to
0$ when $x\to\infty$.

Now define the measure
$$d\mu(s)=\frac{-u'(1/s)}{s^3}1_{[0,1]}(s)ds.$$ An immediate
calculation gives, for $0<\eps\le 1$,
$\int_\eps^1sd\mu(s)=u(1)-u(1/\eps)$, which is bounded, so that
$\int_0^1sd\mu(s)<\infty$. So we can define the Bernstein function
associated to $\mu$: $F(t)=\int_0^\infty(1-e^{-ts})d\mu(s)$. Then,
for all $t\ge 1$, using the inequality $1-e^{-ts}\geq (1-e^{-1})ts$
on $[0,1/t]$:

$$F(t)\ge \int_0^{1/t}(1-e^{-ts})d\mu(s)$$
$$\ge (1-e^{-1})t\int_0^{1/t}sd\mu(s)$$
$$= (1-e^{-1})t\int_0^{1/t}\frac{-u'(1/s)ds}{s^2}$$
$$= (1-e^{-1})t\,u(t)$$
$$\ge (1-e^{-1})w(t)\qquad\text{for large }t.$$

The Bernstein function $x\mapsto (1-e^{-1})^{-1}F(x)$ satisfies our
purposes, as it is easy to see that it is sublinear.\ep

\medskip

An example of an application of Lemma \ref{largesublin} is the
following result.

\begin{Prop}Let $G$ be a locally compact, compactly generated
group having a 1-cocycle of linear growth (e.g. $G=\Z^n$ or $\R^n$
for $n\ge 1$). Then, for every function $f:G\to\R_+$ with sublinear
growth, there exists on $G$ a sublinear 1-cocycle $b$ such that
$\|b(g)\|\succeq f(g)$.
\end{Prop}

\bpr Let $b'$ denote a 1-cocycle with linear growth, and write
$|g|=\|b'(g)\|$, so that $g\mapsto|g|$ is equivalent to the word
length, and its square is conditionally negative definite on $G$.

By hypothesis, $f(g)\prec |g|$. Define $w:\R_{+}\rightarrow \R_{+}$
by
$$w(x)= \sup\{f(h): |h|\leq x\}.$$
Then $w$ is sublinear on $\R_{+}$, and so is the function $x\mapsto
w(x^{1/2})^2$. By Lemma \ref{largesublin}, we find a sublinear
Bernstein function $F$ such that $F(x)\ge w(x^{1/2})^2$ for
large~$x$. Using Lemma \ref{lem:Bochner_Schoenberg}, the function
$g\mapsto F(|g|^2)$ is conditionally negative definite on $G$;
moreover $F(|g|^2)^{1/2}\prec|g|$, and $F(|g|^2)^{1/2}\geq f(g)$ for
$g\in G$ with $|g|$ large enough.\penalty10000\ep

\section{Cocycles with non-slow growth}\label{Sec:lin_growth}

\subsection{Amenability}

Here is a generalization of a result by Guentner and Kaminker
\cite[§5]{GK} who proved it in the case of discrete groups.

\begin{Thm}
Let $G$ be a locally compact group, and $S$ a symmetric, compact
generating subset. Suppose that $G$ admits a 1-cocycle $b$ with
compression $\rho(r)\succ r^{1/2}$. Then $G$ is
amenable.\label{thm:GK_lc}
\end{Thm}

\begin{Cor}
If a locally compact, compactly generated group admits a linear
1-cocycle, then it is amenable.\ep
\end{Cor}

\textbf{Proof of Theorem \ref{thm:GK_lc}} For $t>0$, define
$f_t(g)=e^{-t\|b(g)\|^2}$. By Schoenberg's Theorem \cite[Appendix
C]{BHV}, $f_t$ is positive definite. We claim that $f_t$ is square
summable. Denote $S_n=\{g\in G:\,|g|_S=n\}$, and fix a left Haar
measure $\mu$ on~$G$. There exists $a<\infty$ such that $\mu(S_n)\le
e^{an}$ for all $n$. Since $\rho(r)\succ r^{1/2}$, there exists
$n_0$ such that, for all $n\ge n_0$, and all $g\in S_n$,
$2t\|b(g)\|^2\ge (a+1)n$. Then, for all $n\ge n_0$,
$$\int_{S_n}f_t(g)^2d\mu(g)=\int_{S_n}e^{-2t\|b(g)\|^2}d\mu(g)$$
$$\le \int_{S_n}e^{-(a+1)n}d\mu(g)\le
\mu(S_n)e^{-(a+1)n}\le e^{-n}.$$ Therefore, the sequence
$(\int_{S_n}f_t(g)^2d\mu(g))$ is summable, so that $f_t$ is
square-summable.

By \cite[Théorème 13.8.6]{Dix}, it follows that there exists a
positive definite, square-summable function $\varphi_t$ on $G$ such
that $f_t=\varphi_t\ast\varphi_t$, where $\ast$ denotes convolution.
In other words, $f_t=\langle \lambda(g) \varphi_t,\varphi_t\rangle$,
where $\lambda$ denotes the left regular representation of $G$ on
$L^2(G)$. Note that $f_t$ converges to 1, uniformly on compact
subsets, when~$t\to 0$. We conclude that $(\varphi_t)_{t>0}$
provides a net of almost invariant unit vectors for the regular
representation of $G$ in $L^2(G)$, so that $G$ is amenable.\ep

\subsection{Cocycles with linear growth}\label{Subs:HFD}

Let us recall a property introduced by Shalom in \cite{Shalom04}: a
group has Property $H_{FD}$ if every unitary representation such
that $\overline{H^1}(\Gamma,\pi)\neq 0$ has a finite-dimensional
subrepresentation.

Here are a few useful results about Property $H_{FD}$.

\begin{enumerate}
    \item  [1)] Property $H_{FD}$ is a quasi-isometry invariant among \textit{discrete}
     amenable groups (Shalom, \cite[Theorem 4.3.3]{Shalom04}).

    \item  [2)] A finitely generated amenable group with Property $H_{FD}$ has a finite index subgroup
with infinite abelianization \cite[Theorem 4.3.1]{Shalom04}.

    \item  [3)] A connected amenable Lie group has Property
    $H_{FD}$ (F. Martin, \cite[Theorem 3.3]{Mar}).
    A polycyclic group has Property $H_{FD}$ \cite[Theorem
    5.1.4]{Shalom04}. Both results rely on a deep result due to
    Delorme \cite[Corollaire~V.2]{Del}: connected solvable Lie groups have
    Property $H_{FD}$.

    \item  [4)] The semidirect product $\Z[1/mn]\rtimes_{m/n}\Z$, and the wreath
    product $F\wr\Z$, where $F$ is any finite group, have Property
    $H_{FD}$ \cite[Theorems 5.2.1 and 5.3.1]{Shalom04}.
    Semidirect products $\left(\R\oplus\bigoplus_{p\in
    S}\Q_p\right)\rtimes_{\frac{m}{n}}\Z$ or
    $\left(\bigoplus_{p\in S}\Q_p\right)\rtimes_{\frac{m}{n}}\Z$, with
    $m,n$ co-prime integers, and $S$ a finite set of prime numbers
    dividing $mn$ also have Property $H_{FD}$ \cite[Proof of Theorem 5.3.1]{Shalom04}.

    \item  [5)] The wreath product $\Z\wr\Z$ does \textit{not} have
    Property $H_{FD}$ \cite[Theorem 5.4.1]{Shalom04}.

\end{enumerate}

We prove:

\begin{Thm}\label{FHDcoclineaires}

Let $G$ be a locally compact, compactly generated group with
Property $H_{FD}$. Suppose that $G$ admits a unitary representation
$\pi$ such that $\textnormal{lin}(G,\pi)$ is nonempty. Then, $G$ has
a compact normal subgroup $K$ such that $G/K$ is isomorphic to some
closed subgroup of $\textnormal{Isom}(\R^n)$. In particular,

\begin{itemize}\item $G$ is quasi-isometric to $\R^m$ for some unique
$m$, \item If $G$ is discrete, then $G$ is virtually
abelian.\end{itemize}
\end{Thm}

{\bf Proof:} Let $(\HH,\pi)$ be a unitary representation of $G$ and
suppose that there exists $b\in \text{lin}(G,\pi)$. Replacing $G$ by
$G/K$ for some compact normal subgroup if necessary, we can suppose
by \cite[Theorem 3.7]{Com} that $G$ is separable, and thus we can
also assume that $\HH$ is separable.

As $G$ has Property $H_{FD}$, $\HH$ splits into a direct sum
$\HH=\HH'\oplus(\bigoplus_{n\in \mathbf{N}} \HH_n)$ where $\HH_n$
are finite dimensional subrepresentations and where $\HH'$ is a
subrepresentation with trivial reduced cohomology. By Proposition
\ref{coboundaries}, and since $b$ has linear growth, its orthogonal
projection on $\oplus_{n\in \mathbf{N}} \HH_n$ still has linear
growth, so we can assume that $\HH=\oplus_{n\in \mathbf{N}} \HH_n$.
Now, let $b_n$ be the projection of $b$ on $\oplus_{k\leq n} \HH_k$.
Then $b_n\rightarrow b$ uniformly on compact subsets. So, as
$\text{lin}(G,\pi)$ is open and $b\in \text{lin}(G,\pi)$, there
exists $n$ such that $b_n\in \text{lin}(G,\pi)$. Hence $b_n$ defines
a proper morphism $G\to\Isom(\mathcal{H}_n)$; denote by $H$ its
image.

If $G$ is discrete, by Bieberbach's Theorem (see for instance
\cite{Buser}), this implies that $G$ has a homomorphism with finite
kernel onto a virtually abelian group, hence is itself virtually
abelian.

In general, by Corollary \ref{cor:bieb_loc_cpt}, $G$ acts properly
and cocompactly on some Euclidean space $\R^n$, hence is
quasi-isometric to $\R^n$. \hfill$\square$

\subsection{Uniform embeddings into Hilbert
spaces}\label{subs:unif_emb}

Let $G$ be a locally compact group, and $|\cdot|_S$ the length
function with respect to a compact symmetric generating subset $S$.
For an arbitrary map $f$ of $G$ to a Hilbert space $\HH$, define its
\textit{dilation} as:
$$\delta(x)=\sup\{\|f(g)-f(h)\|:\;|g^{-1}h|_S\le x\}\,\in\R_+\cup\{\infty\},$$
and its \textit{compression} as:
$$\rho(x)=\inf\{\|f(g)-f(h)\|:\;|g^{-1}h|_S\ge x\}\,\in\R_+\cup\{\infty\},$$

We call $f$ a \textit{uniform map} if $\delta(x)<\infty$ for all
$x\in\R_+$ (by an easy standard argument, this implies that $\delta$
has at most linear growth). The map $f$ is called a \textit{uniform
embedding} if, in addition, $\rho(x)\to\infty$ when $x\to\infty$. It
is called a quasi-isometric embedding if, in addition, it has
compression with linear growth, i.e. $\rho(r)\succeq r$.

\medskip

The following result, which was pointed out to us by M. Gromov (who
provided a proof in the discrete case), is very useful.

\begin{Prop}
Let $G$ be a locally compact, compactly generated, amenable group.
Let $f$ be a uniform map of $G$ into a Hilbert space, and $\rho$ its
compression, $\delta$ its dilation. Then $G$ admits a 1-cocycle with
compression $\ge\rho-a$ and dilation $\le\delta +a$, for some
constant $a\ge 0$. If $G$ is discrete, we can take
$a=0$.\label{prop:rem_gromov}
\end{Prop}

\bpr Let $m$ be a mean on $G$, that is, a continuous, linear map on
$L^\infty(G)$ such that $m(1)=1$, $m(f)\ge 0$ whenever $f\ge 0$
locally almost everywhere. Since $G$ is amenable, we choose $m$ to
be right invariant, i.e. $m(f\cdot g)=m(f)$ for all $g\in G$ and
$f\in L^\infty(G)$, where $(f\cdot g)(h)$ is by definition equal
to~$f(hg^{-1})$.

For $g,h\in G$, set $\Psi(g,h)=\|f(g)-f(h)\|^2$. By assumption,
$$\rho(|g^{-1}h|_S)\le
\Psi(g,h)^{1/2}\le\delta(|g^{-1}h|_S),\quad\forall g,h\in G.$$

By Lemma \ref{lem:prol_hilb} in the appendix, there exists a
uniformly continuous function $f'$ at bounded distance from~$f$ (if
$G$ is discrete we do not need Lemma \ref{lem:prol_hilb} since it
suffices to take $f'=f$). Write $\Psi'(g,h)=\|f'(g)-f'(h)\|^2$. Then
$\Psi^{1/2} - (\Psi')^{1/2}$ is bounded.

Now set $u_{g_1,g_2}(h)=\Psi'(hg_1,hg_2)$ for $g_{1}, g_{2}, h\in
G$. The upper bound by $\delta$ and the uniform continuity of $f'$
imply that the mapping $(g_1,g_2)\mapsto u_{g_1,g_2}$ is a
continuous function from $G\times G$ to $L^\infty(G)$, so that the
function $\Psi_m(g_1,g_2)=m(u_{g_1,g_2})$ is continuous on $G\times
G$.

For $g_1,...,g_n\in G$ and $\lambda_1,...,\lambda_n\in\R$ with
$\sum_{i=1}^n \lambda_i=0$, we have
$$\sum_{i,j}\lambda_i\lambda_j\Psi_m(g_i,g_j)=\sum_{i,j}\lambda_i\lambda_j
m(u_{g_i,g_j})=m(\sum_{i,j}\lambda_i\lambda_j u_{g_i,g_j})
\leq 0$$ because $m$ is positive and $\sum_{i,j}\lambda_i\lambda_j
u_{g_i,g_j}$ is a non-positive function on $G$ (since $\Psi'$ is a
conditionally negative definite kernel). So $\Psi_m$ is a
conditionally negative definite kernel.

Since $m$ is right $G$-invariant, it follows that $\Psi_m$ is
$G$-invariant, so that we can write
$\Psi_m(g_1,g_2)=\psi(g_1^{-1}g_2)$ for some continuous,
conditionally negative definite function $\psi$ on $G$. Let $b$ be
the corresponding 1-cocycle. The estimates on $\psi$, and thus on
$\|b\|$, follow from the positivity of~$m$.\ep

\begin{Cor}
If a locally compact, compactly generated amenable group $G$
quasi-isometrically embeds into a Hilbert space, then it admits a
1-cocycle with linear growth.\ep\label{cor:rem_gromov_linear}
\end{Cor}

From Corollary \ref{cor:rem_gromov_linear} and Theorem
\ref{FHDcoclineaires}, we deduce immediately:

\begin{Cor} Let $G$ be a locally compact, compactly generated amenable group with
property $H_{FD}$. The group $G$ admits a quasi-isometric embedding
into a Hilbert space if and only if $G$ acts properly on a Euclidean
space. In particular, if $G$ is discrete, this means that it is
virtually abelian.
\end{Cor}

Combining this corollary with Shalom's results mentioned in
§\ref{Subs:HFD}, we immediately obtain Theorem
\ref{thm:intro_noQIemb} in the introduction.

\medskip

\textbf{Proof of Corollary \ref{cor:rigidite_virt_ab}}. We must
prove that if a finitely generated group $\Gamma$ is
quasi-isomorphic to $\Z^n$, then it has a finite index subgroup
isomorphic to $\Z^n$. This result is known as a consequence of
Gromov's polynomial growth Theorem (see e.g. \cite[Th\'eor\`eme
1.17]{GhHar}, with a sketch of proof on p. 13); it has been given a
new proof in \cite{Shalom04}. As in \cite{Shalom04}, the first step
is the fact that, since Property $H_{FD}$ is a quasi-isometric
invariant of amenable groups, $\Gamma$ has Property $H_{FD}$. Now,
being quasi-isometric to $\Z^n$, $\Gamma$ quasi-isometrically embeds
into a Hilbert space, hence is virtually isomorphic to $\Z^m$ for
some $m$, by Theorem \ref{thm:intro_noQIemb}.

Finally, it is well-known that $\Z^m$ and $\Z^n$ being
quasi-isometric implies $m=n$. For instance, it suffices to observe
that the degree of growth of $\Z^n$ is~$n$.\ep

\medskip

\textbf{Proof of Corollary \ref{cor:tree_not_QI_embed_Hilb}}. It is
enough to show that the regular tree of degree 3 does not embed
quasi-isometrically into a Hilbert space. But such a tree is
quasi-isometric to $G=\Q_2\rtimes_{2}\Z$, since this group acts
cocompactly and properly on the Bass-Serre tree of
$\text{SL}_2(\Q_2)$. On the other hand, $\Q_2\rtimes_{2}\Z$ has no
non-trivial compact normal subgroup (indeed, such a subgroup
would be contained in $\Q_2$, in which the non-trivial conjugacy
classes of $G$ are unbounded), so it does not act properly
cocompactly on a Euclidean space. Since it has Property $H_{FD}$, by
Theorem \ref{thm:intro_noQIemb}, it does not quasi-isometrically
embed into a Hilbert space. \hfill$\square$

\begin{Rem}
Proposition \ref{prop:rem_gromov} is specific to amenable groups.
For instance, if $\Gamma$ is a free group on $n\ge 2$ generators,
then it has uniform embeddings with compression $\ge |g|^a$ for
arbitrary $a<1$ \cite[§6]{GK}, while it has no 1-cocycle with
compression $\succ |g|^{1/2}$ since it is non-amenable
(\cite[§5]{GK} or Theorem \ref{thm:GK_lc}).
\end{Rem}

\appendix

\section{Maps into Hilbert spaces}

Let $G$ be a locally compact group and $\HH$ a Hilbert space. Let
$f$ be a map: $G\to\HH$ (not necessarily continuous). We call $f$ a
\textit{uniform map} if, for every compact subset $K\subset G$, we
have $\sup_{g\in G,k\in K}\|f(kg)-f(g)\|<\infty$. If $G$ is
compactly generated, this coincides with the definition given in
§\ref{subs:unif_emb}.

\begin{Lem}
Let $f:G\to\HH$ be a uniform map\footnote{It follows from the proof
that we can take any Banach space instead of $\HH$.}. Then there
exists $\tilde{f}:G\to\HH$ such that:\begin{itemize}\item
$\tilde{f}$ is at bounded distance from $f$, and\item $\tilde{f}$ is
uniformly continuous on $G$.\end{itemize}\label{lem:prol_hilb}
\end{Lem}

\bpr Fix an open, relatively compact, symmetric neighbourhood $V$ of
$1$ in $G$. Consider a closed, discrete subset $X\subset G$ such
that \begin{itemize}\item[(1)] $\bigcup_{x\in X}xV=G$, and
\item[(2)] for all $x,y\in X$, if $x^{-1}y\in V$, then
$x=y$.\end{itemize}

The existence of such a subset $X$ is immediate from Zorn's Lemma.

Fix a function $\phi:G\to\R_+$, continuous with compact support,
such that $\phi\le 1$, and (3): $\phi\equiv 1$ on $V$. Fix a
symmetric, compact subset $W$ containing the support of~$\phi$.

Set $\Phi(g)=\sum_{x\in X}\phi(x^{-1}g)$ and observe that
$\Phi(g)\geq 1$ as a consequence of (1) and (3). Define

$$\tilde{f}(g)=\frac{1}{\Phi(g)}\sum_{x\in
X}\phi(x^{-1}g)f(x).$$

Let us first check that $\tilde{f}$ is at bounded distance from $f$.
For all $g\in G$,

$$\tilde{f}(g)-f(g)=\sum_{x\in
X}\frac{\phi(x^{-1}g)}{\Phi(g)}(f(x)-f(g)).$$

Since $f$ is a uniform map, there exists $M<\infty$ such that for
all $g,h\in G$, $h^{-1}g\in W$ implies $\|f(h)-f(g)\|\le M$. It
follows that, for all $g\in G$, we have $\|\tilde{f}(g)-f(g)\|\le
M$.

Now let us show that $f$ is uniformly continuous. Consider a
neighbourhood $V_0$ of $1$ in $G$ such that $V_0^2\subset V$. It
immediately follows that, for every $g\in G$, the set $X\cap gV_0$
contains at most one element. On the other hand, by compactness,
there exist $g_1,\dots,g_n$ such that
$W\subset\bigcup_{i=1}^ng_iV_0$. It follows that, for all $g\in G$,
the set $gW\cap X$ has cardinality at most $n$.

Write $u_\phi(g)=\sup_{h\in G}|\phi(h)-\phi(hg)|$. Since $\phi$ is
uniformly continuous, $u_\phi(g)\to 0$ when $g\to 1$.

Then $|\phi(x^{-1}g)-\phi(x^{-1}h)|\le u_\phi(g^{-1}h)$ and, for all
$g,h\in G$, $\Phi(g)-\Phi(h)=\sum_{x\in
X}(\phi(x^{-1}g)-\phi(x^{-1}h))\le 2n\,u_\phi(g^{-1}h)$, since the
only nonzero terms are those for $x\in (gW\cap X)\cup(hW\cap X)$.
Accordingly, $\Phi$ is uniformly continuous. Since $\Phi\ge 1$, it
follows that $1/\Phi$ is also uniformly continuous; let us define
$u_{1/\Phi}$ as we have defined $u_\phi$.

For $g,h\in G$,

$$\left|\frac{\phi(x^{-1}g)}{\Phi(g)}-\frac{\phi(x^{-1}h)}{\Phi(h)}\right|\le
\left|\frac{\phi(x^{-1}g)}{\Phi(g)}-\frac{\phi(x^{-1}h)}{\Phi(g)}\right|+
\left|\frac{\phi(x^{-1}h)}{\Phi(g)}-\frac{\phi(x^{-1}h)}{\Phi(h)}\right|$$
$$\le\frac{|\phi(x^{-1}g)-\phi(x^{-1}h)|}{\Phi(g)}+
|\phi(x^{-1}h)|\left|\frac{1}{\Phi(g)}-\frac{1}{\Phi(h)}\right|$$
$$\le u_{\phi}(g^{-1}h)+u_{1/\Phi}(g^{-1}h).$$

Therefore, fixing some $x_0\in X$,

$$\|\tilde{f}(g)-\tilde{f}(h)\|=\sum_{x\in
X}\left(\frac{\phi(x^{-1}g)}{\Phi(g)}-\frac{\phi(x^{-1}h)}{\Phi(h)}\right)(f(x)-f(x_{0}))$$
$$\le \sum_{x\in
X}\left|\frac{\phi(x^{-1}g)}{\Phi(g)}-\frac{\phi(x^{-1}h)}{\Phi(h)}\right|\|f(x)-f(x_{0})\|$$
$$\le \sum_{x\in
(gW\cap X)\cup(hW\cap
X)}\left|\frac{\phi(x^{-1}g)}{\Phi(g)}-\frac{\phi(x^{-1}h)}{\Phi(h)}\right|\|f(x)-f(x_{0})\|$$
$$\le (u_{\phi}(g^{-1}h)+u_{1/\Phi}(g^{-1}h))\sum_{x\in
(gW\cap X)\cup(hW\cap X)}\|f(x)-f(x_{0})\|.$$

Since $f$ is a uniform map, there exists $M'<\infty$ such that
$h^{-1}g\in V^2W$ implies $\|f(h)-f(g)\|\le M'$ for all~$g,h\in G$.
Now fix $x_0$ so that $g\in x_0V$, and suppose $g^{-1}h\in V$. If
$x\in gW\cup hW$, then it follows that $\|f(x)-f(x_0)\|\le M'$.
Accordingly, whenever $g^{-1}h\in V$,

$$\|\tilde{f}(g)-\tilde{f}(h)\|\le 2n(u_{\phi}(g^{-1}h)+u_{1/\Phi}(g^{-1}h))M',$$

so that $\tilde{f}$ is uniformly continuous.\ep

\section{Actions on Euclidean spaces}

\begin{Prop}
Let $G$ be a closed subgroup of
$\textnormal{E}_n(\R)=\textnormal{Isom}(\R^n)$. The following are
equivalent:\begin{itemize} \item[(i)] $G$ is cocompact in
$\textnormal{E}_n(\R)$.\item[(ii)] $G$ acts cocompactly on
$\R^n$.\item[(iii)] $G$ does not preserve any proper affine subspace
of $\R^n$.\end{itemize}\label{prop:cocompact_euclidien}
\end{Prop}
\bpr (i)$\Leftrightarrow$(ii)$\Rightarrow$(iii) are trivial.

Let us show (iii)$\Rightarrow$(ii). We use some results of Guivarc'h
on the structure of closed (not necessarily connected) subgroups of
amenable connected Lie groups. By \cite[Théorème IV.3 and Lemma
IV.1]{Guiv}, $G$ has a characteristic closed cocompact solvable
subgroup $R$. Then $R$ has a characteristic subgroup of finite index
$N$ which maps to a torus of $\text{O}_n(\R)$ through the natural
projection $\textnormal{E}_n(\R)\to\textnormal{O}_n(\R)$.

\textit{First case}: $G$ does not contain any nontrivial
translation. Then $N$ is abelian. If $g\in N$, let $d_g$ denote its
displacement length: $d_g=\inf\{\|g.v-v\|:\,v\in\R^n\}$. The subset
$A_g=\{v\in\R^n:\,\|gv-v\|=d_g\}$ is a (nonempty) affine subspace of
$\R^n$, and is $N$-stable since $N$ is abelian. Also note that if
$W$ is any $g$-stable affine subspace, then\footnote{Let $v\in A_g$
such that $d(v,W)=d(A_g,W)$. Let $p$ denote the projection on $W$;
since $W$ is $g$-stable, $p$ commutes with $g$. Since
$d(v,pv)=d(gv,gpv)\le d(x,y)$ for all $x\in [v,gv],y\in[pv,gpv]$, we
easily obtain that $v$, $gv$, $pv$, $gpv$ form a rectangle, so that
$d(pv,gpv)=d(v,gv)$ and thus $pv\in A_g$ by definition of $A_g$.}
$W\cap A_g\neq\emptyset$. It easily follows that finite
intersections of subspaces of the form $A_g$, for $g\in N$, are
nonempty, and a dimension argument immediately yields that
$A=\bigcap_{g\in N}A_g\neq\emptyset$. This is a $G$-invariant affine
subspace, hence is, by assumption, all of $\R^n$. Therefore, every
element in $N$ is a translation, so that $N=\{1\}$ and thus $G$ is
compact. This implies that $G$ has a fixed point, so that the
assumption implies $n=0$ (i.e. leads to a contradiction if $n\ge
1$).

\textit{General case}. Argue by induction on the dimension $n$.
Suppose that $n\ge 1$. Let $T_G$ be the subgroup of translations in
$G$. Let $W$ be the linear subspace generated by $T_G$. Since $T_G$
is closed, it acts cocompactly on $W$. Moreover, by the first case,
$W$ has positive dimension. Note that the linear action of $G$
clearly preserves~$W$.

Now look at the action of $G$ on the affine space $\R^n/W$. It does
not preserve any proper affine subspace, hence is cocompact by the
induction hypothesis. Since the action of $T_G$ on $W$ is also
cocompact, it follows that the action of $G$ on $\R^n$ is also
cocompact.

\begin{Cor}Let $G$ be a locally compact group. Suppose that $G$
has a proper isometric action on a Euclidean space. Then it has a
proper cocompact isometric action on a Euclidean
space.\label{cor:bieb_loc_cpt}
\end{Cor}
\bpr Let $G$ act on a Euclidean space by isometries. Let $V$ be a
$G$-invariant affine subspace of minimal dimension. Then the action
of $G$ on $V$ is clearly proper, and is cocompact by Proposition
\ref{prop:cocompact_euclidien}.\ep

\bigskip
\footnotesize
\noindent Yves de Cornulier\\
IRMAR, Campus de Beaulieu, \\
F-35042 Rennes Cedex, France\\
E-mail: \url{decornul@clipper.ens.fr}\\
\medskip

\noindent Romain Tessera\\
Department of mathematics, Vanderbilt University,\\ Stevenson
Center, Nashville, TN 37240,
USA,\\
 E-mail: \url{tessera@clipper.ens.fr}
\medskip

\noindent Alain Valette\\
Institut de Mathématiques - Université de Neuchâtel\\
Rue Emile Argand 11, CH-2007 Neuchâtel, Switzerland\\
E-mail: \url{alain.valette@unine.ch}

\end{document}